\documentclass[11pt,reqno]{amsart}
\usepackage{graphicx}
\usepackage{mathtools}
\usepackage{ulem}
\usepackage{color}

\usepackage{amsmath}
\usepackage{amssymb}
\usepackage{graphicx}
\usepackage{accents}
\usepackage{pstricks,pst-node,pst-coil,pst-plot}
\usepackage{mathtools}
\usepackage{ulem}
\usepackage{mathtools}
\usepackage{ulem}
\usepackage{graphicx}
\usepackage{color}
\usepackage{graphicx}
\usepackage{color}
\usepackage{dcolumn}
\usepackage{graphicx}
\usepackage{color}
\usepackage{dcolumn}
\usepackage{bm}
\usepackage{mathrsfs}
\usepackage{graphicx}
\usepackage{color}
\usepackage{dcolumn}
\usepackage{bm}
\usepackage{slashed}
\usepackage{epstopdf}
\usepackage{slashed}
\usepackage{bm}
\usepackage{slashed}
\usepackage{dcolumn}
\usepackage{bm}
\usepackage{slashed}
\usepackage{color}
\usepackage{color}

\newtheorem{thm}{Theorem}
\newcommand{\be}{\begin{equation}}
\newcommand{\ee}{\end{equation}}

\newtheorem{remark}{Remark}
\newtheorem{lem}{Lemma}
\newtheorem{prop}{Proposition}
\newtheorem{definition}{Definition}

\newcommand{\beq}{\begin{equation}}
\newcommand{\eeq}{\end{equation}}
\newcommand{\bee}{\begin{equation*}}
\newcommand{\eee}{\end{equation*}}
\def\d{\mathrm{d}}

\let\<\langle
\let\>\rangle

\newcommand{\definedas}{\mathrel{\raise.095ex\hbox{\rm :}\mkern-5.2mu=}}

\begin{document}
%%%%%%%%%%%%%%%%%%%%%%%%%%%%%%%%%%%%%%%%%%%%%%%%%%%%%%%%%%%%%%%%%

%%%%%%%%%%%%%%%%%%%%%%%%%%%%%%%%%%%%%%%%%%%%%%%%%%%%%%%%%%%%%%%%%
%\begin{center}
%  \framebox{\framebox{
%      \vbox{This is project {\red \Project}\\
%        Current version {\blue\Version}, from {\blue\Datum}, most
%        recent changes by {\blue\Person}.}  }}
%\end{center}
%%%%%%%%%%%%%%%%%%%%%%%%%%%%%%%%%%%%%%%%%%%%%%%%%%%%%%%%%%%%%%%%%

\title[Stability of smooth periodic travelling waves in DGH equation]{Stability of smooth periodic travelling waves in the Dullin-Gottwald-Holm equation}

\author{Xiaokai He}
\address[Xiaokai He]{School of Mathematics and Statistics, Hunan First Normal University, Changsha 410205, China}
\email{sjyhexiaokai@hnfnu.edu.cn}

\author{Aiyong Chen* }
\address[Aiyong Chen]{School of Mathematics and Statistics, Hunan First Normal University, Changsha 410205, China}
\email{aiyongchen@163.com}

\author{Gengrong Zhang }
\address[Gengrong Zhang]{School of Mathematics and Statistics, Hunan First Normal University, Changsha 410205, China}
\email{18152853519@163.com}

\thanks{$^*$Corresponding author.}

\begin{abstract}
The existence of smooth periodic traveling solutions in the Dullin-Gottwald-Holm (DGH) equation and the monotonicity of the period function are clarified. By introducing two suitable parameters, we show the existence of periodic travelling solutions of DGH equation in a concise way. The monotonicity of period function with respect to different variables are proved by using Chicone's criterion and the method developed by Geyer and Villadelprat.
The problem of the spectral stability of smooth periodic waves in the DGH equation is discussed. Within the functional-analytic framework, we obtain a criterion for the spectral stability of smooth periodic traveling waves
in DGH equation. In addition, we show the smooth periodic travelling solutions are orbitally stable under certain conditions.

\end{abstract}

\subjclass[2010]{35Q35,35P30,37K45}
%
% Latex changes the classification to 1991 ?!  it should be 2010!

\date{\today}

%\date{May 29, 2023}

\keywords{Dullin-Gottwald-Holm equation, smooth periodic travelling waves, spectral stability, orbital stability.}

\maketitle

%\tableofcontents

%%%%%%%%%%%%%%%%%%%%%%%%%%%%%%%%%%%%%%%%%%%%%%%%%%%%%%%%%%%%%%%%%%%%%%%%%%%%%%%%
%%%%%%%%%%%%%%%%%%%%%%%%%%%%%%%%%%%%%%%%%%%%%%%%%%%%%%%%%%%%%%%%%%%%%%%%%%%%%%%%

\section{Introduction}\label{S1}
In 1993, Camassa and Holm derived the so called Camassa-Holm (CH) equation
\begin{equation}
u_t+2\omega u_x-u_{txx}+3uu_x=2u_xu_{xx}+uu_{xxx}
\end{equation}
as a model of the propagation of unidirectional shallow water waves\cite{CH1993}. CH equation is integrable and possesses bi-Hamiltonian structure, Lax pairs and infinite hierachy conservation laws\cite{Fokas1981}. Since the work of Camassa and Holm, many researchers pay attentions to construct  nonlocally evolutive equations of the type (1). Particularly, we address the Dullin-Gottwald-Holm (DGH) equation
\begin{equation}\label{DGH}
m_t+2\omega u_x+um_x+2mu_x+\gamma u_{xxx}=0,
\end{equation}
where $m(x,t)=u(x,t)-\alpha^2\partial_x^2u(x,t)$ is a momentum variable, the constant $c_0=\sqrt{gh}$ (where $c_0=2\omega$) is the linear wave speed for undisturbed water at rest at spatial infinity with $g$ the gravitational constant and $h$ the mean fluid depth, the constants $\alpha^2$ and $\frac{\gamma}{c_0}$ are squares of length scales\cite{Dullin2001}.  The DGH equation was derived by using asymptotic expansions directly in the Hamiltonian for Euler's equations in the shallow water regime. Equation (\ref{DGH}) combines the linear dispersion of the Korteweg-de Vries (KdV) equation with the nonlinear and nonlocal dispersion of the CH equation, and preserves integrability. When $\alpha=1,\gamma=0$, the DGH equation reduces to the CH equation. When $\alpha=0,\gamma=1$, the DGH equation becomes to the KdV equation.

Many researches on the DGH equation had been done in the past few decades. Tian, Gui and Liu \cite{Tian2005} studied the well-posedness of the Cauchy problem and the scattering problem for the DGH equation. Yin and Tian \cite{Yin2010} studied  the peaked solitons for the periodic DGH equation on the line.  They  demonstrated that the solutions of the Cauchy
problem for DGH equation converge to those of the corresponding periodic
Camassa-Holm equation as the linear dispersive parameter converges to zero. Zhou et al\cite{Zhou2013} obtained a two-peakon solution to the DGH equation explicitly
by direct computation, and then discussed the peakon-antipeakon interaction. The classification of bounded travelling wave solutions for the DGH equation is obtained in Ref.\cite{Silva2019}.

Stability of travelling wave solutions is an important aspect in the study of shallow water wave equations. For  peakon solutions, Constantin and Strauss \cite{Constantin2000} proved the peakons of the CH equation were  orbitally stable.  The instability results of $H^1$-stable peakons in the CH equation were obtained in \cite{Natali2020-1}. A variational approach for proving the orbital stability of the peakons was introduced by Constantin and Molient \cite{Constantin2001}.   Moreover, the orbital stability of the single
peakons for the Degasperis-Procesi (DP) equation was proved by Lin and Liu \cite{Lin}. They developed the
approach due to Constantin and Strauss \cite{Constantin2000} in a delicate way.  The approach in \cite{Constantin2000}
was extended in \cite{Liu1} to prove the orbital stability of the peakons for the Novikov
equation. By constructing a function only depending on
three important conservative laws, Yin and Tian proved that the shapes of peakons for DGH equation are stable
under small perturbations\cite{Yin2010}.  
In 2022, Lafortune and Pelinovsky studied the spectral instability of peakons in the b-family of the CH equations \cite{Lafortune2022-1}. For soliton solutions,
 Constantin and Strauss \cite{Constantin2002} proved the orbital stability of smooth solitary wave solutions of the CH equation by using Grillakis-Shatah-Strauss method \cite{Grillakis}. Li, Liu and Wu \cite{Liliuwu1,Liliuwu2} proved the spectral stability and orbital stability of smooth solitary waves for the DP equation. For the periodic travelling wave solutions,  by using the variational method, the orbital stability of peaked periodic waves in $H^1_{\rm per}$ was obtained in 2004\cite{Lenells2004-1,Lenells2004-2}. Here we use $H^s_{\rm per}$ to denote the Sobolev space $H^s_{\rm per}(\mathbb{T}_L)$ with $\mathbb{T}_L\equiv[0,L]$ the periodic domain of length $L>0.$ The stability of smooth periodic solution to CH equation 
was first proven by Lenells\cite{Lenells2005}.

In the functional-analytic framework, stability of smooth periodic traveling waves is a difficult problem. Only very recently, Geyer, Martins, Natali and  Pelinovsky successfully solved the open problem of the spectral stability of smooth periodic waves in the CH equation\cite{Geyer2021}.  Based on the key observation that the periodic waves of the CH equation can be characterized by two different Hamiltonian structure, they showed the nonstandard formulation is suitable for the study of the spectral and orbital stability of the smooth periodic solutions in CH equation.  By constructing two linear operators $\mathcal{L}$ and $\mathcal{K}$ and utilizing the properties of the two operators, they obtained the spectral stability criterion of periodic traveling waves in CH equation. Furthermore, they proved the periodic waves with profile $\phi
\in H^{\infty}_{\rm per}$ is orbitally stable in certain parameter regions.  In 2022, Geyer and Pelinovsky studied the stability of smooth periodic travelling waves in DP equation\cite{Geyer2022}.     Inspired by these new developments, we focus on the stability of smooth periodic traveling waves in DGH equation (2) in this paper.

In terms of $u$, the DGH equation can be rewritten as
\be\label{DGH-u}
(1-\alpha^2\partial_x^2)u_t+2\omega u_x+3uu_x+\gamma u_{xxx}=\alpha^2(2u_xu_{xx}+uu_{xxx}),
\ee
i.e.,
\be
u_t=J[2\omega u+\frac{3}{2}u^2+\gamma u_{xx}-\frac{1}{2}\alpha^2u_x^2-\alpha^2uu_{xx}],
\ee
where
\be\label{J}
J:=-(1-\alpha^2\partial_x^2)^{-1}\partial_x
\ee
is a well-defined operator from $H_{\rm per}^s$ to $H_{\rm per}^{s+1}$ for $s\in\mathbb{R}.$

Let $u:\mathbb{R}\times\mathbb{R}\rightarrow\mathbb{R}$ with $L-$periodic at space variable be a smooth solution to the DGH equation. It is well known that the DGH equation (\ref{DGH-u}) admits the following conserved quantities\cite{Tian2005}
\be\label{Mu}
M(u)=\int_0^L u\d x,
\ee
\be\label{Eu}
E(u)=\frac{1}{2}\int_0^L(u^2+\alpha^2u_x^2)\d x,
\ee
and
\be\label{Fu}
F(u)=\frac{1}{2}\int_0^L(u^3+\alpha^2uu_x^2+2\omega u^2-\gamma u_x^2)\d x.
\ee
The DGH equation has the following Hamiltonian structure 
\be
\frac{\d u}{\d t}=J\frac{\delta F}{\delta u}.
\ee

In this work, we consider the smooth travelling solution to the DGH equation. The smooth travelling solution with speed $c$ and profile $\phi$ takes the form
\be\label{phi}
u(x,t)=\phi(x-ct):=\phi(\xi), \ \ \xi=x-ct.
\ee
Plugging (\ref{phi}) into (\ref{DGH-u}), one can get
\be\label{phi1}
(-\alpha^2c-\gamma+\alpha^2\phi)\phi'''
+2\alpha^2\phi'\phi''+(c-2\omega-3\phi)\phi'=0.
\ee
Direct integration of (\ref{phi1}) in $x$ yields
\be\label{intg1}
(-\alpha^2c-\gamma+\alpha^2\phi)\phi''+\frac{1}{2}\alpha^2\phi'^2
+(c-2\omega-\frac{3}{2}\phi)\phi=b,
\ee
where $b$ is an integration constant.
Let
\be
C_1=c+\frac{\gamma}{\alpha^2},\ C_2=2\omega+\frac{\gamma}{\alpha^2},
\ee
then we have
\be\label{intg1new}
\alpha^2(\phi-C_1)\phi''+\frac{1}{2}\alpha^2\phi'^2+(C_1-C_2-\frac{3}{2}\phi)\phi=b.
\ee

\begin{remark}
The introduction of parameters $C_1$ and $C_2$ is technically important for the study of the smooth periodic travelling waves in DGH equation, since it makes the similarity between the smooth periodic travelling waves in DGH equation and that in CH equation more clear. Using two parameters $C_1$ and $C_2$ rather than the original four parameters $\alpha,\omega, \gamma$ and $c$, enable us to study the existence of periodic travelling solutions in DGH equation in a concise way.
\end{remark}

Direct calculation shows the linearized operator for the second-order equation (\ref{intg1new}) is given by
\be\label{mathcalL}
\mathcal{L}=-\alpha^2\partial_x(C_1-\phi)\partial_x+(C_1-C_2-3\phi+\alpha^2\phi'').
\ee

The linearized operator $\mathcal{L}:L_{\rm per}^2\rightarrow L_{\rm per}^2$ is a self-dual, unbounded operator in $L_{\rm per}^2$ equipped with the standard inner product
\be
\langle f,g\rangle:=\int_0^Lf(x)g(x)\d x.
\ee

From (\ref{intg1new}), one can construct a planar dynamical system as
\be\label{planar1}
\left\{
\begin{aligned}
\frac{\d\phi}{\d \xi}&=X,\\
\frac{\d X}{\d \xi}&=\frac{b-\frac{1}{2}\alpha^2X^2-(C_1-C_2-\frac{3}{2}\phi)\phi}{\alpha^2(\phi-C_1)},
\end{aligned}
\right.
\ee
 which has a first integral
\be\label{intg2new}
\begin{split}
\alpha^2(\phi-C_1)\phi'^2
+(C_1-C_2)\phi^2-\phi^3-2b\phi+C_1(2b+C_1C_2)=C_3,
\end{split}
\ee
with $C_3$  an integration constant.

Combining (\ref{intg1new}) and (\ref{intg2new}), it can be obtained
\be
-2(\phi-C_1)^2(\alpha^2\phi''-\phi)+C_2(\phi-2C_1)\phi=C_3-C_1^2C_2.
\ee
i.e.,
\be\label{intg-phipp}
\alpha^2\phi''=\phi-\frac{1}{2(\phi-C_1)^2}\bigg[
C_3-C_1^2C_2-C_2\phi^2+2C_1C_2\phi\bigg].
\ee

The first result of this paper is related to the existence of smooth periodic travelling wave solutions with profile $\phi$ satisfying (\ref{intg2new}).

\begin{thm}\label{thm1}
For fixed $C_1$ and $C_2$ which satisfy $2C_1+C_2>0$, smooth periodic travelling wave solutions of (\ref{intg2new}) with profile $\phi\in H^{\infty}_{\rm per}$ exist in an open, simply connected region on the $(C_3,b)$ plane enclosed by the following three boundaries:

(1) $C_3=0$ and $b\in(-\frac{C_1^2}{2}-C_1C_2,\frac{1}{8}(C_2^2-4C_1C_2))$,

(2) $C_3=C_3^+(b)$ and $b\in(\frac{1}{8}(C_2^2-4C_1C_2),\frac{1}{6}(C_1-C_2)^2)$,

(3) $C_3=C_3^-(b)$ and $b\in(-\frac{1}{2}C_1^2-C_1C_2,\frac{1}{6}(C_1-C_2)^2)$,$\hspace{32mm}$
where $C_3^+(b)$ and $C_3^-(b)$ are $C^1$ function of $b$,see Fig.\ref{region} for illustration.

\begin{figure}[htp]
\begin{center}
\includegraphics[width=0.80\textwidth]{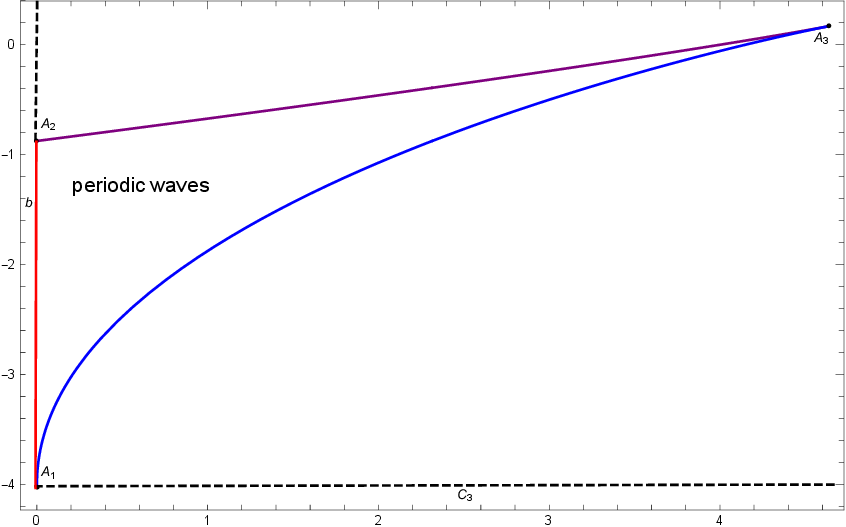}
\caption{\label{region} The existence region of smooth periodic travelling wave solutions on the parameter plane $(C_3,b)$ for $C_1=2$ and $C_2=1$. The points $A_1,A_2$ and $A_3$  are $(0,-\frac{1}{2}C_1^2-C_1C_2)$, $(0,\frac{1}{8}(C_2^2-4C_1C_2))$ and $(\frac{1}{27}(2C_1+C_2)^3,\frac{1}{6}(C_1-C_2)^2)$, respectively. }
\end{center}
\end{figure}

\end{thm}

Secondly,  the monotonicity of the period function for smooth periodic travelling wave solutions in DHG equation is studied. Let
\be
C_{3,critical}=\frac{1}{27}(2C_1+C_2)^3.
\ee

Compared  to  the period function in CH equation\cite{Geyer2021}, the expression of the  period function for smooth travelling solutions in DGH equation is  more complicated due to the presence of the additional parameter $C_2$, which makes the study on the properties of period function in DGH equation more difficult.  Considering this fact, we divide the parameter region $2C_1+C_2>0$ into four subregions in $C_1-C_2$ plane and prove the following theorem by using the monotonicity criterion obtained in \cite{Chicone1987}.

\begin{thm}\label{thm2}
For fixed $C_3\in (0,C_{3,critical})$ and $C_1,C_2$ with $2C_1+C_2>0$ , the period function $\mathfrak{L}(C_1,C_2,C_3,b)$  for smooth periodic travelling wave solutions in DHG equation  is strictly increasing with respect to $b$.
\end{thm}

Moreover, based on the method developed by Geyer and Villadelprat \cite{Geyer2015}, we obtain
\begin{thm} \label{thm3}
For fixed $C_1$ and $C_2$ which satisfy $2C_1+C_2>0$, the period function $\mathfrak{L}(C_1,C_2,C_3,b)$ for smooth periodic travelling wave solutions in DHG equation

(i) is monotonically increasing in $C_3$ if $b\in (-\frac{1}{2}C_1^2-C_1C_2),b_1)$

(ii) has a single maximum point in $C_3$ if $b\in (b_1,\frac{1}{8}(C_2^2-4C_1C_2))$

(iii) is monotonically decreasing in $C_3$ if $b\in(\frac{1}{8}(C_2^2-4C_1C_2),\frac{1}{6}(C_1-C_2)^2)$

where $$b_1=(-1+\frac{\sqrt{6}}{3})C_1^2
+(-\frac{3}{2}+\frac{\sqrt{6}}{3})C_1C_2+(-\frac{1}{8}+\frac{\sqrt{6}}{12})C_2^2.
$$
\end{thm}

Recall that the smooth periodic travelling wave $\phi\in H^{\infty}_{\rm per}$ is said to be spectrally stable in the evolution  problem (\ref{DGH-u}) if the spectrum of $J\mathcal{L}$ in $L^2_{\rm per}$ is located on the imaginary axis, where $J$ and $\mathcal{L}$ are given by (\ref{J}) and (\ref{mathcalL}). We consider the energy criterion for the spectral stability of the periodic waves in the DGH equation (\ref{DGH-u}) and obtain

\begin{thm}\label{thm4}
For fixed $C_1, C_2$ which satisfy $2C_1+C_2>0$ and a fixed period $L>0$, there exists a $C^1$ mapping $C_3\mapsto\mathfrak{B}_L(C_3)$ for $C_3\in(0,C_3^L)$ with some $C_3^L\in (0,C_{3,critical})$ and a $C^1$ mapping $C_3\mapsto\phi=\Phi_L(\cdot,C_3)\in H^{\infty}_{\rm per}$ for smooth $L-$periodic solutions along the curve $b=\mathfrak{B}_L(C_3)$. Let
\be
\begin{split}
\mathcal{M}_L(C_3):=M(\Phi_L(\cdot,C_3)),\\
\mathcal{F}_L(C_3):=F(\Phi_L(\cdot,C_3)),\\
\end{split}
\ee
where $M(u)$ and $F(u)$ are given by (\ref{Mu}) and (\ref{Fu}). The $L-$periodic wave with profile $\phi=\Phi_L(\cdot,C_3)$ is spectrally stable if the mapping
\be
C_3\mapsto\frac{\mathcal{F}_L(C_3)}{\mathcal{M}_L(C_3)^3}
\ee
is strictly decreasing in $C_3$ and $\partial_{C_3}\mathfrak{L}(C_1,C_2,C_3,b)<0$.
\end{thm}

\begin{remark}
When $2C_1+C_2<0$, similar results to Theorems 1,2,3 and 4  can still be obtained.
\end{remark}

The rest of this paper is organized as follows. We devote Section 2 to a discussion of the existence of smooth travelling solutions in DGH equation
and give the proof of Theorem 1. The monotonicity of the smooth periodic travelling waves is then studied in Section 3. The proof of Theorem 2 and Theorem 3, and the limit behaviour of the period function are given in this Section. In Section 4, the spectral and orbital stability of the smooth periodic travelling waves are investigated. Based on the relation between the eigenvalues of the linear operator $\mathcal{L}$ and the monotonicity of the period function, we prove Theorem 4 in this section. Furthermore, the orbital stability of the periodic travelling waves in DGH equation is discussed.

\section{Existence of smooth travelling wave solutions}

In this section, we discuss the existence of smooth periodic travelling wave solutions in DGH equation and provide the proof of Theorem \ref{thm1}.

Firstly, we can rewrite (\ref{intg2new}) as
\be\label{newton1}
b=\frac{1}{2}\alpha^2\phi'^2+U(\phi),
\ee
where
\be
U(\phi)=-\frac{1}{2}\phi^2-\frac{1}{2}C_2\phi
-\frac{1}{2}C_1C_2-\frac{C_3}{2(\phi-C_1)}.
\ee
Let \be\label{z}
z=\frac{\xi}{\alpha},
\ee then (\ref{newton1}) can be rewritten as
\be\label{newton2}
b=\frac{1}{2}(\frac{\d\phi}{\d z})^2+U(\phi),
\ee
which can be seen as a system with total energy $b$ of Newton's particle of unit mass with coordinate $\phi$ in ``time" z with the potential energy $U(\phi).$

The existence of smooth periodic solution is determined by the profile of the potential function $U(\phi)$. Direct calculation shows
\be
\begin{split}
\frac{\d U}{\d \phi}%=&-\phi-\frac{1}{2}C_2+\frac{2A+C_1(2b+C_1C_2)}{2(\phi-C_1)^2}\\
=&\frac{1}{2(\phi-C_1)^2}\bigg[-2(\phi-C_1)^2(\phi+\frac{1}{2}C_2)+C_3\bigg].
\end{split}
\ee

Let
\be
\begin{split}
f(\phi)\equiv 2(\phi-C_1)^2(\phi+\frac{1}{2}C_2),
\end{split}
\ee
then
\be
f'(\phi)=2(\phi-C_1)(3\phi+C_2-C_1)
\ee
which implies the critical points of $f(\phi)$ are
\be
\phi_{1}^{*}=C_1,\ \phi_2^{*}=\frac{C_1-C_2}{3}.
\ee
Moreover, we have
\be
f''(\phi)=12\phi-8C_1+2C_2,
\ee
and hence
\be
\begin{split}
f''(\phi_1^*)&=4C_1+2C_2,\ \
 f''(\phi_2^*)=-(4C_1+2C_2).
\end{split}
\ee

\begin{figure}[htp]
\begin{center}
\includegraphics[width=0.50\textwidth]{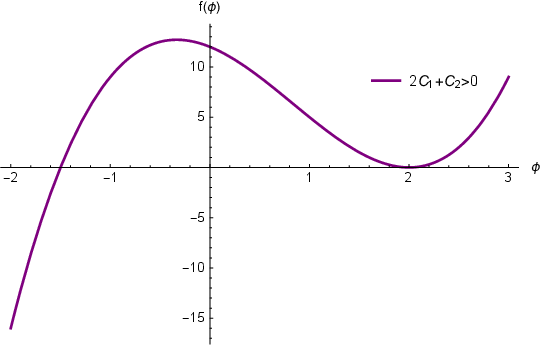}
\caption{\label{f1} Graph of $f(\phi)$ with $2C_1+C_2>0$. }
\end{center}
\end{figure}

When $2C_1+C_2>0$, the local maximum of $f(\phi)$ occurs at $\phi=\frac{C_1-C_2}{3}$ (See Fig.\ref{f1}), and
\be
f_{max}=f(\frac{C_1-C_2}{3})=\frac{1}{27}(2C_1+C_2)^3.
\ee

Let
\be
C_{3,critical}\equiv \frac{1}{27}(2C_1+C_2)^3.
\ee

\begin{figure}[htp]
\begin{center}
\includegraphics[width=0.50\textwidth]{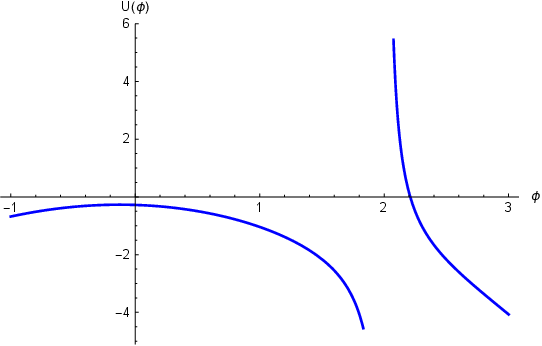}
\caption{\label{U1-2} Graph of $U(\phi)$ with $2C_1+C_2>0$ and $C_3\in(-\infty,0)$ . }
\end{center}
\end{figure}

\begin{figure}[htp]
\begin{center}
\includegraphics[width=0.9\textwidth]{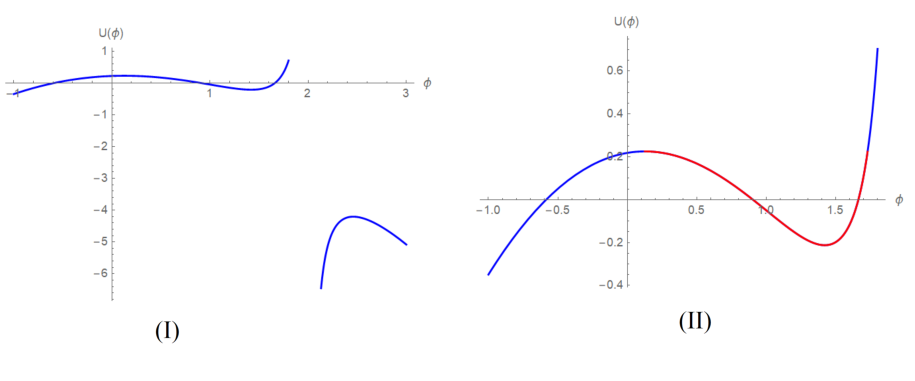}
\caption{\label{U1} Graph of $U(\phi)$ with $2C_1+C_2>0$ and
 $C_3\in(0,C_{3,critical})$. Fig.(II) is just the left branch of Fig.(I). The red portion corresponds to the existence region of smooth periodic travelling waves. }
\end{center}
\end{figure}

\begin{figure}[htp]
\begin{center}
\includegraphics[width=0.50\textwidth]{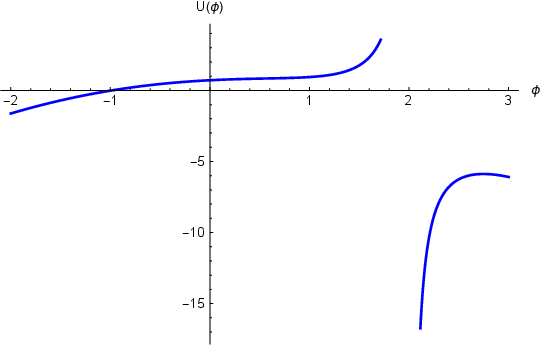}
\caption{\label{U1-3} Graph of $U(\phi)$ with $2C_1+C_2>0$ and $C_3\in(C_{3,critical},\infty)$ . }
\end{center}
\end{figure}

 For $C_3\in(-\infty,0)\cup(C_{3,critical},+\infty)$, $f(\phi)=C_3$ has only one solution which implies there exists only one critical point of $U(\phi)$.
For $C_3\in(0,C_{3,critical})$, $f(\phi)=C_3$ has only three solutions which implies there exists three critical points of $U(\phi)$, two are local maximum and one is local mimimum. Furthermore, $\phi=C_1$ is the pole singularity of $U(\phi)$ if $C_3\neq 0$. See Figs.\ref{U1-2},\ref{U1},\ref{U1-3} for illustration.


Follows from the dynamics of system (\ref{newton1}) or (\ref{newton2}), we can obtain all smooth mapping $x\mapsto\phi(x)$ for
$C_3\in(-\infty,0)\cup(C_{3,critical},+\infty)$ are unbounded.
For $C_3\in (0,C_{3,critical})$, we label the critical points of $U(\phi)$ as $\phi_1<\phi_2<\phi_3$, where $\phi_1,\phi_2,\phi_3$ are also the roots of the equation $f(\phi)=C_3$ which satisfying the ordering
\be\label{ordering-1}
-\frac{1}{2}C_2<\phi_1<\frac{C_1-C_2}{3}<\phi_2<C_1<\phi_3.
\ee

The local minimum of $U(\phi)$ at $\phi_2$ gives the center of smooth travelling equation at $(\phi_2,0)$ enclosed by the homoclinic orbit connecting the saddle $(\phi_1,0)$.  More precisely, the phase portrait can be discussed as follows.

From (\ref{planar1}), it can be obtained
\be\label{planar1}
\left\{
\begin{aligned}
\frac{\d\phi}{\d z}&=Y,\\
\frac{\d Y}{\d z}&=\frac{b-\frac{1}{2}Y^2-(C_1-C_2-\frac{3}{2}\phi)\phi}{\phi-C_1},
\end{aligned}
\right.
\ee
which has a first integral
\be\label{first1}
(\phi-C_1)\dot{\phi}^2+(C_1-C_2)\phi^2-\phi^3-2b\phi+C_1(2b+C_1C_2)=C_3,
\ee
where $\dot\phi=\frac{\d\phi}{\d z}.$ Let
\be
\hat{f}(\phi)\equiv b-(C_1-C_2-\frac{3}{2}\phi)\phi=\frac{3}{2}\phi^2-(C_1-C_2)\phi+b.
\ee
The roots of $\hat{f}(\phi)=0$ are
\be
\phi_{1}=\frac{C_1-C_2-\sqrt{(C_1-C_2)^2-6b}}{3}
\ee
and
\be
\phi_{2}=\frac{C_1-C_2+\sqrt{(C_1-C_2)^2-6b}}{3}.
\ee
Inserting $\phi=C_1$ into (\ref{first1}) yields
\be
C_3(C_1)=0.
\ee
Inserting $\phi=\phi_1,\dot\phi=0$ into (\ref{first1}) yields
\be
C_3(\phi_1)=\frac{1}{27}\big(2C_1+C_2+\sqrt{(C_1-C_2)^2-6b} \big)g(C_1,C_2,b),
\ee
where
\be
g(C_1,C_2,b)=12b+2C_1^2-C_2^2+8C_1C_2-(2C_1+C_2)\sqrt{(C_1-C_2)^2-6b}.
\ee
By the theory of dynamical system, the kinds of phase portrait  can be classified according to the sign of $C_3(\phi_1)-C_3(C_1)$, or equivalently, the sign of $g(C_1,C_2,b)$.

We show the phase portrait on the phase plane $(\phi,\dot{\phi})$ for

$ (i)\  C_1=3,C_2=1.02, b=-1$, in which case $g=7.32894>0,$

$(ii)\  C_1=2.01,C_2=0.03,b=-1$, in which case $g=-16.1944<0,$

\begin{figure}[htp]
\begin{center}
\includegraphics[width=0.50\textwidth]{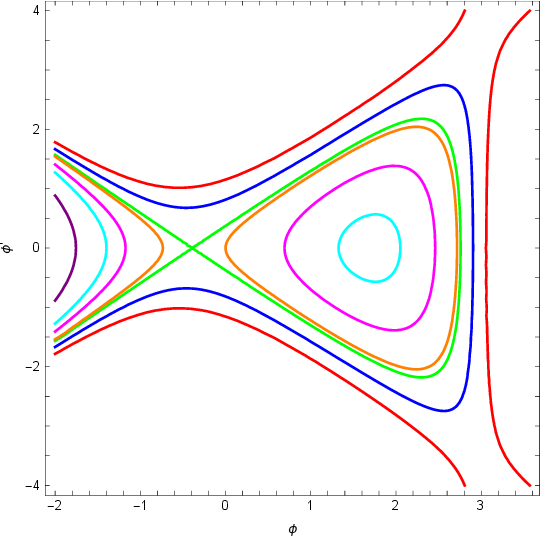}
\caption{\label{DGH1} Phase portrait constructed from the level curves of the first order invariant (\ref{first1}) on the phase plane for  $C_1=3,C_2=1.02, b=-1$. }
\end{center}
\end{figure}

\begin{figure}[htp]
\begin{center}
\includegraphics[width=0.50\textwidth]{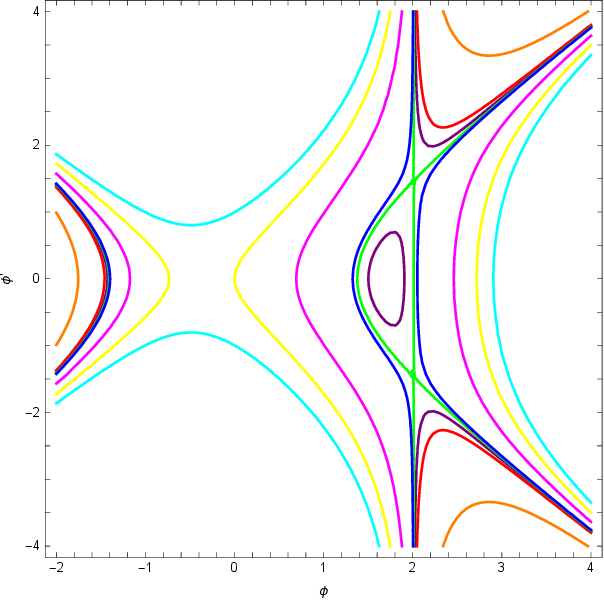}
\caption{\label{DGH2} Phase portrait constructed from the level curves of the first order invariant (\ref{first1}) on the phase plane for $C_1=2.01,C_2=0.03,b=-1$. }
\end{center}
\end{figure}

\begin{figure}[htp]
\begin{center}
\includegraphics[width=0.50\textwidth]{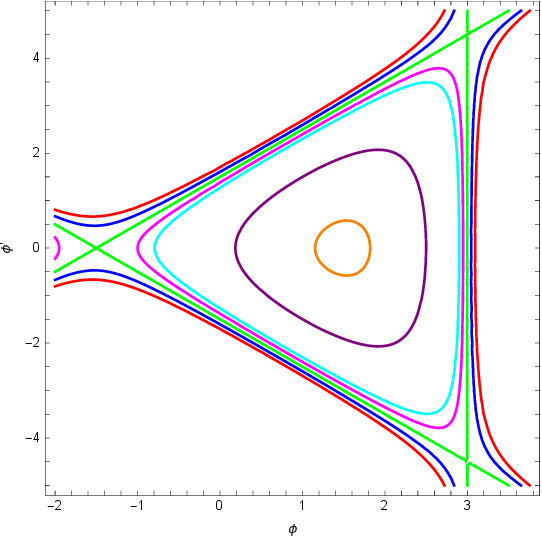}
\caption{\label{DGH3} Phase portrait constructed from the level curves of the first order invariant (\ref{first1}) on the phase plane for $C_1=3,C_2=3,b=-27/8$. }
\end{center}
\end{figure}

$(iii)\  C_1=3,C_2=3,b=-27/8$, in which case $g=0,$

\noindent in Fig.\ref{DGH1} , Fig.\ref{DGH2} and Fig.\ref{DGH3}, from which one can see the period annulus around the center.

 The smooth periodic solutions for fixed $C_3\in (0,C_{3,critical})$ and $C_1,C_2$ with $2C_1+C_2>0$ are parameterized by the parameter $b\in(b_-,b_+)$, where
\be
b_-=U(\phi_2),\ \  b_+=U(\phi_1).
\ee
Especially, for $C_3=0$, we can get
\be
b_-=-\frac{C_1^2}{2}-C_1C_2,\ b_+=\frac{1}{8}(C_2^2-4C_1C_2).
\ee
For $C_3=C_{3,critical}$, we have
\be
b_-=b_+=\frac{1}{6}(C_1-C_2)^2.
\ee

Noting that the first order invariant (\ref{newton1}) is smooth with respect to parameters $C_1,C_2,C_3$ and $b$, the periodic orbits inside the period annulus are also smooth with respect to these parameters.

For fixed $C_1, C_2$ which satisfy $2C_1+C_2>0$, the boundary $b=b_1(C_3)$ corresponds to the center, hence $\phi(x)=\phi_2$ is constant in $x$. Along the curve $b=b_-(C_3)$, $C_3$ and $b$ can be parameterized by $\phi_2$ as
\be\label{bc3}
\begin{split}
b=&(C_1-C_2-\frac{3}{2}\phi_2)\phi_2,\\
C_3=&(C_1-\phi_2)^2(C_2+2\phi_2),
\end{split}
\ee
which results from (\ref{intg1new}) and (\ref{intg2new}) by using $\phi=\phi_2$ is constant. By (\ref{bc3}) we can obtain
\be
\frac{\d b}{\d \phi_2}=C_1-C_2-3\phi_2,\ \frac{\d C_3}{\d\phi_2}=2(C_1-C_2-3\phi_2)(C_1-\phi_2),
\ee
which implies
\be
\frac{\d C_3}{\d b}=2(C_1-\phi_2)>0.
\ee
Hence the mapping $C_3\mapsto b_-(C_3)$ is $C^1$, invertible and monotonous increasing from $(C_3,b)=(0,-\frac{1}{2}C_1^2-C_1C_2)$ to $(C_3,b)=(\frac{1}{27}(2C_1+C_2)^3,\frac{1}{6}(C_1-C_2)^2)$.
Similarly, one can show the mapping $C_3\mapsto b_+(C_3)$ is $C^1$, invertible and monotonous increasing from
 $(C_3,b)=(0,\frac{1}{8}(C_2^2-4C_1C_2))$ to $(C_3,b)=(\frac{1}{27}(2C_1+C_2)^3,\frac{1}{6}(C_1-C_2)^2)$.

Denote the inverse of $b=b_-(C_3)$ and $b=b_+(C_3)$ by $C_3=C_3^-(b)$ and $C_3=C_3^+(b)$, then the assertions in theorem 1 holds.

\section{Monotonicity of the period function for smooth periodic travelling waves }

Now we study the monotonicity of the period function for smooth periodic travelling wave solutions in DGH equation. Theorem 2 and 3 will be proved in this section.

Let us now define the period function for the smooth periodic waves of theorem \ref{thm1}. For a given smooth periodic travelling wave solution, by (\ref{newton2}), the period function $\mathfrak{L}(C_1,C_2,b,C_3)$ is defined by the following integral
\be\label{periodfunction-2}
 \mathfrak{L}(C_1,C_2,b,C_3):= 2\int_{\phi_-}^{\phi_+}\sqrt{\frac{1}{2[b-U(\phi)]}}\d\phi,
\ee
where $\phi_-$ and $\phi_+$ are the turning points of the Newton system (\ref{newton2}) which satisfy the ordering
\be\label{order2}
-\frac{1}{2}C_2<\phi_1<\phi_-<\phi_2<\phi_+<C_1<\phi_3.
\ee

Recall that we have introduced the variable $z$ by (\ref{z}), then
(\ref{intg-phipp}) becomes
\be\label{phidd-1}
\frac{\d^2\phi}{\d z^2}=\phi-\frac{1}{2(\phi-C_1)^2}\bigg[
C_3-C_1^2C_2-C_2\phi^2+2C_1C_2\phi\bigg].
\ee

\begin{figure}[htp]
\begin{center}
\includegraphics[width=0.55\textwidth]{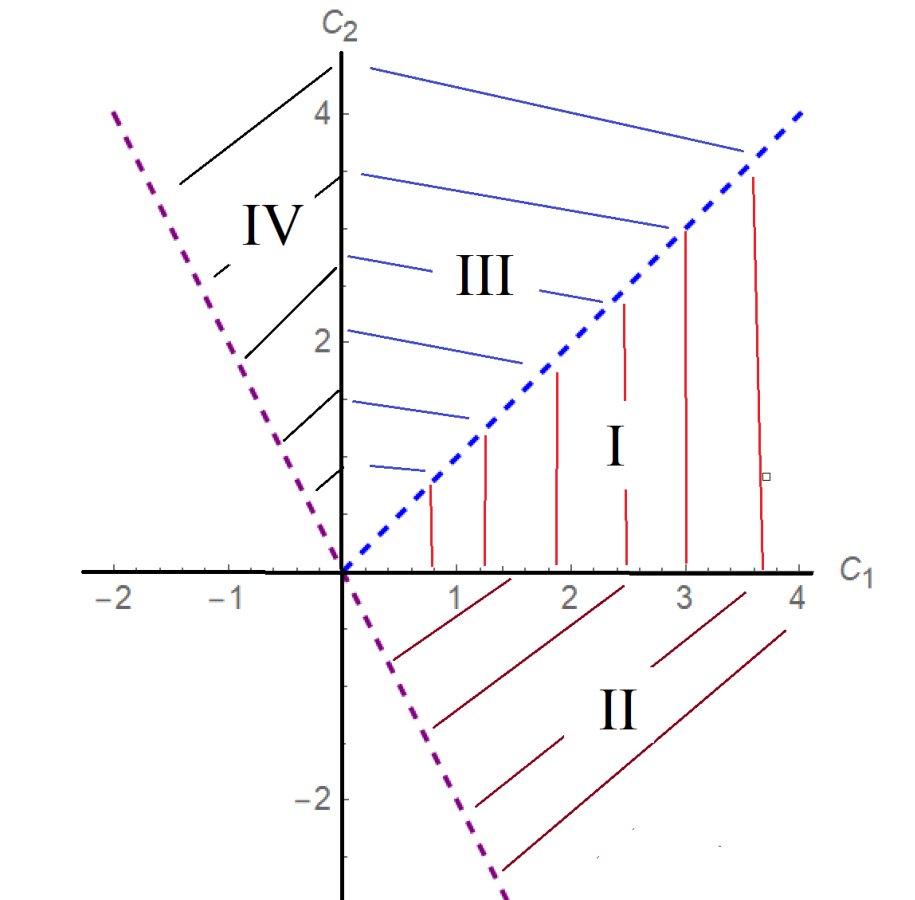}
\caption{\label{c1c2} Four regions in the $C_1-C_2$ plane. }
\end{center}
\end{figure}

Let $Q$ be the second root of $C_3=f(\phi)$ in the ordering (\ref{ordering-1}), i.e., $Q=\phi_2$. In general, $Q$ may be positive, negative or zero. To prove Theorem \ref{thm2}, we divide the region $2C_1+C_2>0$ into 4 parts in the $C_1-C_2$ plane, see Fig.\ref{c1c2}.

\begin{lem}\label{lem1}

For fixed $C_3\in (0,C_{3,critical})$ and $C_1,C_2$ with $C_1>C_2>0$ , the period function $\mathfrak{L}(C_1,C_2,C_3,b)$  defined in (\ref{periodfunction-2}) is strictly increasing with respect to $b$.

\end{lem}

\begin{figure}[htp]
\begin{center}
\includegraphics[width=0.60\textwidth]{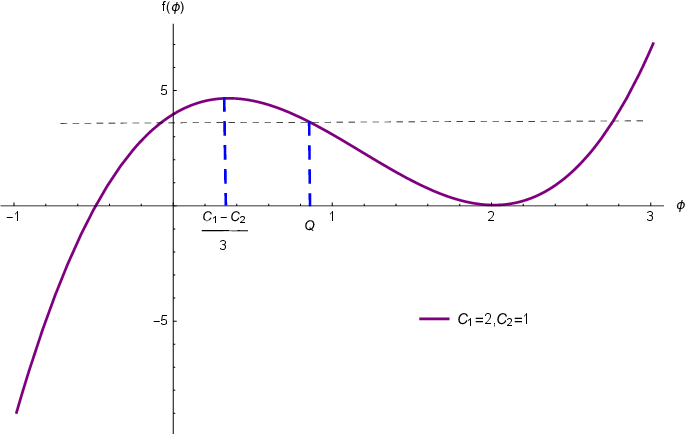}
\caption{\label{f3} Graph of $f(\phi)$ with $C_1=2,C_2=1$. }
\end{center}
\end{figure}

\textit{Proof.} When $C_1>C_2>0$, $Q>\frac{C_1-C_2}{3}$ is positive, see Fig.\ref{f3} for illustration.

By the transformation
\be
x=\frac{\phi-Q}{Q},y=\frac{\dot{\phi}}{Q},
\ee
 the second-order equation (\ref{phidd-1}) can be recast as the planar system
\be\label{newdy}
\begin{split}
\frac{\d x}{\d z}&=y,\\
\frac{\d y}{\d z}&=1+x-\frac{2\eta^2-\beta(x-2\eta)x}{2(x-\eta)^2},
\end{split}
\ee
associated with the Hamiltonian
\be\label{Ham1}
\begin{split}
H(x,y)=&\frac{1}{2}y^2+G(x),\\
G(x)=&-\frac{1}{2}\bigg[x^2
+(\beta+2)x
+\frac{(\beta+2)\eta^2}{x-\eta}+(\beta+2)\eta
\bigg],
\end{split}
\ee
where
\be
\beta=\frac{C_2}{Q},\ \eta=\frac{C_1-Q}{Q}.
\ee

From the relation
$
0<\frac{C_1-C_2}{3}<Q<C_1
$ and the definitions of $\eta$ and $\beta$,
it can be obtained
\be
0<\frac{C_2}{C_1}<\beta<\frac{3C_2}{C_1-C_2},\ 0<\eta<\frac{2C_1+C_2}{C_1-C_2},
\ee
and
\be
\eta-2<\beta<\eta+1.
\ee

\begin{figure}[htp]
\begin{center}
\includegraphics[width=0.80\textwidth]{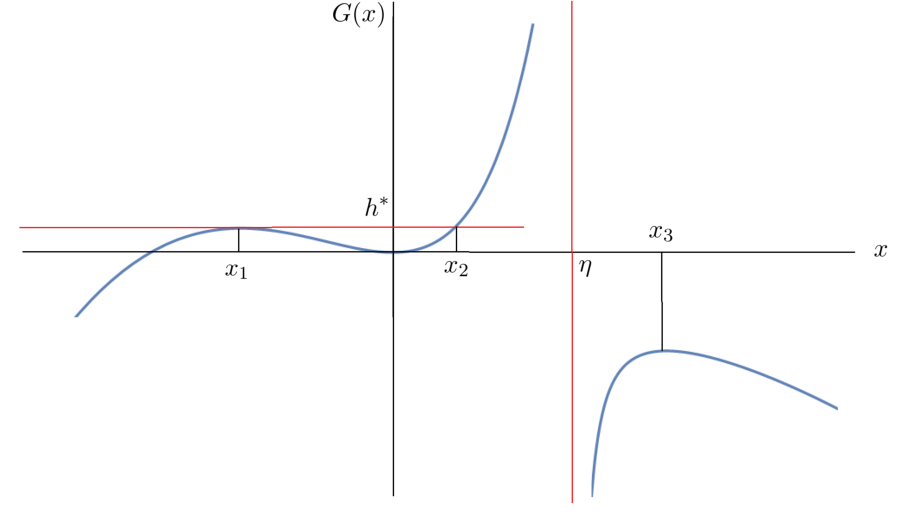}
\caption{\label{Gx} Graph of the potential function $G(x)$. }
\end{center}
\end{figure}

The potential function $G(x)$ is smooth away from the singular line $x=\eta$, see Fig.\ref{Gx}. Moreover, $G(x)$ has a local mimimum at $x=0$ and two maxima at
\be
x_1=\frac{1}{4}\bigg(
4\eta-\beta-2-\sqrt{(\beta+2)(\beta+2+8\eta)}\bigg)
\ee
and
\be
x_3=\frac{1}{4}\bigg(
4\eta-\beta-2+\sqrt{(\beta+2)(\beta+2+8\eta)}\bigg).
\ee

For  dynamical system (\ref{newdy}), the center  at the origin $(x=0,y=0)$ is surrounded by periodic orbits $\gamma_h$, which lie inside the level curves $H(x,y)=h$ with $h\in(0,h^*)$ and
$h^*=G(x_1)$. Denote by $x_2$ the solution of $G(x)=h^*$ which obeys the ordering
\be
x_1<0<x_2<\eta<x_3.
\ee
The period function of the center $(0,0)$ of system (\ref{newdy}) is defined as
\be
\mathbb{L}(h):=\int_{\gamma_h}\frac{\d x}{\d y},\ \ h\in(0,h^*).
\ee
Noting that
\be
h=\frac{b}{Q^2}+\frac{1}{2}-\eta+\beta
\ee
and
\be
\mathfrak{L}(C_1,C_2,C_3,b)= \mathbb{L}(h)
\ee
for fixed $C_3\in(0,C_{3,critical})$ and $C_1>C_2>0.$ Since $Q$ is fixed, then $\partial_b \mathfrak{L}(C_1,C_2,C_3,b)>0$ if and only if $\mathbb{L}'(h)>0.$

To prove  $\mathbb{L}'(h)>0$, we use the monotonicity criterion by Chicone\cite{Chicone1987}.
Let
\be
W(x):=\frac{G(x)}{[G'(x)]^2},
\ee
Chicone showed $\mathbb{L}'(h)>0$ if $W''(x)>0$ for every $x\in (x_1,x_2).$

Direct calculation shows
\be
W(x)=-\frac{2(x-\eta)^3(x+\beta+2\eta)}{[2x^2+(\beta+2-4\eta)x+2\eta(\eta-\beta-2)]^2}
\ee
and
\be
W''(x)=-\frac{12(\beta+2)(x-\eta)R(x)}{[2x^2+(\beta+2-4\eta)x+2\eta(\eta-\beta-2)]^4},
\ee
where
\be
\begin{split}
R(x)=&(4\eta-\beta-2)x^3+\eta(14+7\beta-12\eta)x^2+
3\eta^2(4\eta-6-3\beta)x\\
&+\eta^2(\beta^2+6\eta-4\eta^2+4\beta+3\beta\eta+4).
\end{split}
\ee
Since $x_2<\eta,\beta>0$, we need to show that $R(x)>0$ for every $x\in[x_1,x_2]$, $\ \beta\in
(\frac{C_2}{C_1},\frac{3C_2}{C_1-C_2}),\ \eta\in (0,\frac{2C_1+C_2}{C_1-C_2})$.

When $4\eta-\beta-2=0$, $R(x)$ is simply
\be
R(x)=8\eta^2(2x^2-3\eta x+3\eta^2).
\ee
Due to the fact $\Delta=9\eta^2-24\eta^2=-15\eta^2<0$, $R(x)$ is strictly positive for all $x$.

When $4\eta-\beta-2\neq 0$, the discriminant of  the cubic function $R(x)$ with respect to $x$ is
\be
D=-(\beta+2)^3\eta^4 S(\eta,\beta),
\ee
where
\be
\begin{split}
S(\eta,\beta)=&27\beta^3+2(81+92\eta)\beta^2
+(324+736\eta-240\eta^2)\beta\\
&+8(
27+92\eta-60\eta^2+16\eta^3).
\end{split}
\ee
For every fixed $\eta >0$, $S(\eta,\beta)$ can be seen as a function in $\beta$. Taking the derivative of $S(\eta,\beta)$ with respect to $\beta$ yields
\be\label{pS}
\begin{split}
\frac{\partial S}{\partial\beta}=&
81\beta^2+4(81+92\eta)\beta+(324+736\eta-240\eta^2)\\
=&(\beta-\beta_1)(\beta-\beta_2),
\end{split}
\ee
where
\be
\begin{split}
\beta_1=-\frac{2}{81}(81+92\eta-2\sqrt{3331}\eta),\\
\beta_2=-\frac{2}{81}(81+92\eta+2\sqrt{3331}\eta).
\end{split}
\ee
From  (\ref{pS}), one can obtain
\be
\frac{\partial S}{\partial\beta}\bigg|_{\beta=\beta_1}=
\frac{\partial S}{\partial\beta}\bigg|_{\beta=\beta_2}=0,
\ee
and
\be
\frac{\partial^2 S}{\partial\beta^2}\bigg|_{\beta=\beta_1}=\beta_1-\beta_2>0,\
\frac{\partial^2 S}{\partial\beta^2}\bigg|_{\beta=\beta_2}=\beta_2-\beta_1<0,
\ee
which imply $S(\beta,\eta)$ has one local maximum at $\beta=\beta_2$ and one local minimum at $\beta=\beta_1$  for every fixed $\eta$. Noting that, for any $\eta >0 $,
\be
\begin{split}
S(\beta_1,\eta)=&\frac{128(200854-3331\sqrt{3331})\eta^3}{19683}>0,\\
S(0,\eta)=&8(27+92\eta-60\eta^2+16\eta^3)>0.
\end{split}
\ee
Hence we can conclude that $S(\beta,\eta)$ is strictly positive for $\eta>0,\beta>0.$  It follows the discriminant of  the cubic function $R(x)$ with respect to $x$ is strictly negative for $\ \beta\in
(\frac{3C_2}{C_1-C_2},\frac{C_2}{C_1}), \eta\in (0,\frac{2C_1+C_2}{C_1-C_2})$. Therefore, if the
$4\eta-\beta-2\neq 0$, the cubic polynomical $R(x)$ has exactly one real root, $x_0$.

For $4\eta-\beta-2> 0$ case, it follows the dominant behaviour of $R(x)$ that $R(x)\rightarrow +\infty$ as $x\rightarrow +\infty$ and $R(x)\rightarrow -\infty$ as
$x\rightarrow -\infty$. If
\be
R(x_1)>0,
\ee
then the only real root $x_0$ must locate for $x_0<x_1$. Therefore,
\be
R(x)>0,\ x\in[x_1,x_2].
\ee

To prove $R(x_1)>0$, we note
\be
\begin{split}
R(x_1)=&\frac{1}{16}(\beta+2)(\beta+2+8\eta)\bigg[4+\beta^2+2\sqrt{(\beta+2)(\beta+2+8\eta)}\\
&+\eta\bigg(4-2\sqrt{(\beta+2)(\beta+2+8\eta)}\bigg)\\
&+\beta\bigg(4+2\eta+\sqrt{(\beta+2)(\beta+2+8\eta)}\bigg)\bigg].
\end{split}
\ee
Let
\be
\begin{split}
N(\beta,\eta)=&4+\beta^2+2\sqrt{(\beta+2)(\beta+2+8\eta)}\\
&+\eta\bigg(4-2\sqrt{(\beta+2)(\beta+2+8\eta)}\bigg)\\
&+\beta\bigg(4+2\eta+\sqrt{(\beta+2)(\beta+2+8\eta)}\bigg).
\end{split}
\ee
Then
\be
N(\beta,\eta)=(\beta+2)(\beta+2+2\eta)+(\beta+2-2\eta)\sqrt{(\beta+2)(\beta+2+8\eta)}.
\ee
When $\beta+2-2\eta\geq 0$, it is obviously that $N(\beta,\eta)>0.$

When $\beta+2-2\eta < 0$, then $N(\beta,\eta)>0$ is equivalent to
\be
(\beta+2)(\beta+2+2\eta)>(2\eta-\beta-2)\sqrt{(\beta+2)(\beta+2+8\eta)}
\ee
which holds identically since $\beta+2>\eta$.

Therefore $N(\beta,\eta)>0$ which implies $R(x_1)>0.$

For $4\eta-\beta-2< 0$ case, it follows the dominant behaviour of $R(x)$ that $R(x)\rightarrow +\infty$ as $x\rightarrow -\infty$ and $R(x)\rightarrow -\infty$ as
$x\rightarrow +\infty$. Notice that
\be
R(\eta)=(2+\beta)^2\eta^2>0,
\ee
so the only real root $x_0$ must locate for $x_0>\eta$. Therefore,
\be
R(x)>0,\ x\in[x_1,x_2]
\ee
since $x_2<\eta.$

This completes the proof of Lemma \ref{lem1}. $\ \ \Box$

\begin{lem}\label{lem2}

For fixed $C_3\in (0,C_{3,critical})$ and $C_1,C_2$ which satisfy $C_1>0\geq C_2$ and $2C_1+C_2>0$ , the period function $\mathfrak{L}(C_1,C_2,C_3,b)$  defined in (\ref{periodfunction-2}) is strictly increasing with respect to $b$.

\end{lem}

\begin{lem}\label{lem3}

For fixed $C_3\in (0,C_{3,critical})$ and $C_1,C_2$ which satisfy $C_2\geq C_1>0$ , the period function $\mathfrak{L}(C_1,C_2,C_3,b)$  defined in (\ref{periodfunction-2}) is strictly increasing with respect to $b$.

\end{lem}

\begin{lem}\label{lem4}

For fixed $C_3\in (0,C_{3,critical})$ and $C_1,C_2$ which satisfy $C_2>0, C_1\leq 0$ and $2C_1+C_2>0$ , the period function $\mathfrak{L}(C_1,C_2,C_3,b)$  defined in (\ref{periodfunction-2}) is strictly increasing with respect to $b$.

\end{lem}

By similar arguments used in the proof of Lemma \ref{lem1}, one can verify Lemma \ref{lem2}, Lemma \ref{lem3} and Lemma \ref{lem4}. %See Appendix \ref{app1}, \ref{app2} and \ref{app3} for details.
Hence, the monotonicity result of the period function $\mathfrak{L}(C_1,C_2,C_3,b)$   with respect to $b$ in Theorem \ref{thm2} follows from Lemmas 1,2,3 and 4.

Now we analyze the monotonicity of $\mathfrak{L}(C_1,C_2,C_3,b)$ with respect to $C_3$ and prove Theorem \ref{thm3}. The monotonicity of the period function in $b$ is more subtle than the one with respect to $C_3$. The results can be obtained by using the method developed by  Geyer and Villadelprat in Ref.\cite{Geyer2015}.

Recall that (\ref{intg1new}) can be rewritten as
\be\label{intg1new-2}
(\phi-C_1)\frac{\d^2\phi}{\d z^2}+\frac{1}{2}(\frac{\d\phi}{\d z})^2+(C_1-C_2-\frac{3}{2}\phi)\phi=b.
\ee

Let $w=\phi-C_1,\ v=\frac{\d\phi}{\d z}$, then (\ref{intg1new-2}) yields
\be
\begin{split}
\frac{\d w}{\d z}=&v,\\
\frac{\d v}{\d z}=&-\frac{F'(w)+\frac{1}{2}v^2}{w},
\end{split}
\ee
where
\be
\begin{split}
F(w)=&\alpha w+\beta w^2-\frac{1}{2}w^3,\\
\alpha=&-b-\frac{1}{2}C_1^2-C_1C_2,\\
\beta=&-\frac{1}{2}(2C_1+C_2).
\end{split}
\ee
It is easy to verify that the function
\be
H(w,v)=\frac{1}{2}wv^2+F(w)
\ee
is a first integral of the above planar dynamical system. In terms of $\phi$ and $\frac{\d \phi}{\d z}$, $h$ can be expressed as
\be
\begin{split}
h=&\frac{1}{2}wv^2+F(w)\\
=&\frac{1}{2}\big(C_1-\phi)(2b+C_1C_2+C_2\phi+\phi^2-(\frac{\d\phi}{\d z})^2\big).
\end{split}
\ee
This $h$ is exactly  $\frac{1}{2}C_3$ due to (\ref{intg2new}). Following exactly the same analysis used in Ref.\cite{Geyer2015}, one can show directly that Theorem \ref{thm3} holds, i.e., for fixed $C_1$ and $C_2$ which satisfy $2C_1+C_2>0$, the period function $\mathfrak{L}(C_1,C_2,C_3,b)$ for smooth periodic travelling wave solutions in DHG equation

$\bullet$ is monotonically increasing in $C_3$ if $b\in (-\frac{1}{2}C_1^2-C_1C_2),b_1)$,

$\bullet$ has a single maximum point in $C_3$ if $b\in (b_1,\frac{1}{8}(C_2^2-4C_1C_2))$,

$\bullet$ is monotonically decreasing in $C_3$ if $b\in(\frac{1}{8}(C_2^2-4C_1C_2),\frac{1}{6}(C_1-C_2)^2)$,

where $$b_1=(-1+\frac{\sqrt{6}}{3})C_1^2
+(-\frac{3}{2}+\frac{\sqrt{6}}{3})C_1C_2+(-\frac{1}{8}+\frac{\sqrt{6}}{12})C_2^2.$$

\begin{remark}
When $C_1=c, C_2=0$, the above results about the monotonicity of  $\mathfrak{L}(C_1,C_2,C_3,b)$ with respect to $b$  reduce to Theorem 2.5 in Ref.\cite{Geyer2015} or Remark 1 in Ref.\cite{Geyer2021}.

\end{remark}

In the rest of this section, we discuss the limit behaviour of the period function.

\begin{lem}\label{b-minus}

Fix $C_1,C_2$ with $2C_1+C_2>0$ and $C_3\in(0,C_{3,critical})$. The smooth periodic solutions of Theorem \ref{thm1} transform to the constant solutions as $b\rightarrow b_-(C_3)$. The limiting period function
\be
\mathfrak{L}_-(C_3):=\lim_{b\rightarrow b_-(C_3) }\mathfrak{L}(C_1,C_2,b,C_3)
\ee
is monotonous increasing in $C_3$ along the curve $b=b_-(C_3)$. Moreover,
$\mathfrak{L}_-(C_3)\rightarrow 0$ as $C_3\rightarrow 0$ and $\mathfrak{L}_-(C_3)\rightarrow \infty$ as $C_3\rightarrow C_{3,critical}$.

\end{lem}

\textit{Proof.} From the ordering (\ref{order2}), it can be obtained that the boundary $b=b_-(C_3)$ corresponds to the center $\phi_-=\phi_+=\phi_2$, hence $\phi(x)=\phi_2$ is constant in $x$. By the linearization of the second-order equation (\ref{intg-phipp}) at the center point $(\phi_2,0)$, one can get the period $\mathfrak{L}_-(C_3)$ as
\be\label{omega-1}
\mathfrak{L}_-(C_3)=\frac{2\pi}{\omega}, \ \ \omega=\sqrt{\frac{C_3}{(C_1-\phi_2)^3}-1}.
\ee
As shown previously, along the curve $b=b_-(C_3)$, $C_3$ and $b$ can be parameterized by $\phi_2$ as
\be\label{b-1}
b=(C_1-C_2-\frac{3}{2}\phi_2)\phi_2,
\ee
\be\label{c3-1}
C_3=(C_1-\phi_2)^2(C_2+2\phi_2).
\ee
 Substituting (\ref{c3-1}) into (\ref{omega-1}) yields
\be\label{omega-2}
\omega^2=\frac{3\phi_2-C_1+C_2}{C_1-\phi_2}.
\ee
Combining (\ref{omega-1}) with (\ref{b-1}), it can be deduced
\be\label{b-2}
b=\frac{(C_1-C_2)^2}{6}-\frac{8(2C_1+C_2)^2\pi^4}{3[4\pi^2+3(\mathfrak{L}_-(C_3))^2]^2}
\ee
It follows from (\ref{b-2}) that $\mathfrak{L}_-(C_3)$ increases in $b$ along the curve $b=b_-(C_3)$. Moreover, $\mathfrak{L}_-(C_3)\rightarrow 0$ as $b\rightarrow
-\frac{C_1^2}{2}-C_1C_2$ ( equivalently, $C_3\rightarrow 0$ ) and $\mathfrak{L}_-(C_3)\rightarrow \infty$ as $b\rightarrow \frac{1}{6}(C_1-C_2)^2$ ( equivalently, $C_3\rightarrow C_{3,critical}$). This completes the proof of Lemma \ref{b-minus}.

\begin{lem}\label{b-plus}
Fix $C_1,C_2$ with $2C_1+C_2>0$ and $C_3\in(0,C_{3,critical})$. The smooth periodic solutions of Theorem \ref{thm1} transform to the solitary solutions as $b\rightarrow b_+(C_3)$. The limiting period function
\be
\mathfrak{L}_+(C_3):=\lim_{b\rightarrow b_+(C_3)}\mathfrak{L}(C_1,C_2,b,C_3)=\infty.
\ee

\end{lem}

\textit{Proof.} From the ordering (\ref{order2}), it can be obtained that the boundary $b=b_+(C_3)$ corresponds to $\phi_- =\phi_1$. Hence, $\phi(x)$ is the solitary wave solution which satisfies $\phi(x)\rightarrow\phi_1$
as $x\rightarrow \infty$ so that  $\mathfrak{L}_+(C_3)=\infty$.

\begin{lem}\label{c3=0}

Fix $C_1,C_2$ with $2C_1+C_2>0$ and $C_3\in(0,C_{3,critical})$. The smooth periodic solutions of Theorem \ref{thm1} transform to the peaked periodic solutions as $C_3\rightarrow 0$. The period function
\be
\mathfrak{L}_0(b):=\mathfrak{L}(C_1,C_2,b,0)
\ee
is monotonous increasing in $b$. Moreover,
$\mathfrak{L}_0(b)\rightarrow 0$ as $b\rightarrow -\frac{1}{2}C_1^2-C_1C_2$ and $\mathfrak{L}_0(b)\rightarrow \infty$ as $b\rightarrow\frac{1}{8}(C_2^2-4C_1C_2)$.

\end{lem}

\textit{Proof.} Noting that $\phi_+$ and $\phi_-$ are the turning points of the Newton system (\ref{newton2}), i.e.,
\be\label{phipm}
(C_1-\phi_{\pm})(2b+\phi_{\pm}^2+C_2\phi_{\pm}+C_1C_2)=C_3.
\ee
If $C_3=0$, we can get
\be
\begin{split}
\phi_+=&C_1,\\
\phi_-=&\frac{-C_2+\sqrt{C_2^2-4(2b+C_1C_2)}}{2},
\end{split}
\ee
since $\phi_+$ and $\phi_-$ satisfy the ordering (\ref{order2}).

Moreover, when $C_3=0$, it follows from (\ref{intg-phipp}) that $\phi$ satisfies the ordinary differential equation
\be\label{phi-z1}
\frac{\d^2\phi}{\d z^2}-\phi-\frac{C_2}{2}=0.
\ee
Without loss of generality, one can place the maximum of $\phi$ at $z=0$ and the minimum of $\phi$ at $z=\pm\frac{L}{2}$ which yields
\be\label{phi-z2}
\phi(0)=\phi_+=C_1.
\ee
Equation (\ref{phi-z1}) with condition (\ref{phi-z2}) can be solved explicitly as
\be
\phi(z)=(C_1+\frac{C_2}{2})\frac{\cosh(\frac{L}{2}-|x|)}{\cosh(\frac{L}{2})}-\frac{C_2}{2}.
\ee
This periodic wave is peaked at $z=0$ and smooth at $z=\pm\frac{L}{2}$ with $\frac{\d\phi}{\d z}\big|_{\pm\frac{L}{2}}=0.$ From (\ref{intg1new}), direct calculation shows
\be\label{b-L-2}
b=-\frac{8C_1^2+12C_1C_2+C_2^2+(4C_1-C_2)C_2\cosh L}{16(\cosh\frac{L}{2})^2}
\ee
It follows from (\ref{b-L-2}) that $L=\mathfrak{L}_0(b)$ increasing in $b$ since
\be
\frac{\d b}{\d L}=(2C_1+C_2)^2\frac{(\sinh\frac{L}{2})^4}{(\sinh L)^3}.
\ee
 In addition,
$\mathfrak{L}_0(b)\rightarrow 0$ as $b\rightarrow -\frac{1}{2}C_1^2-C_1C_2$ and $\mathfrak{L}_0(b)\rightarrow \infty$ as $b\rightarrow\frac{1}{8}(C_2^2-4C_1C_2)$. This completes the proof of Lemma \ref{c3=0}.

\section{Stability of the smooth periodic travelling waves}

In this section, we focus on the stability of the smooth periodic travelling wave solutions in DGH equation.

Adding a co-periodic perturbation $v(x,t)$ to the smooth travelling wave $\phi$ propagating with the same fixed speed $c$ in
\be
u(x,t)=\phi(x-ct)+v(x-ct,t).
\ee
At the linear order, one can obtain the perturbation equation from the DGH equation (\ref{DGH})
\be
v_t=-J\mathcal{L}v.
\ee
where $J$ and $\mathcal{L}$ are given by (\ref{J}) and (\ref{mathcalL}) respectively.

\begin{definition}
 The smooth periodic travelling wave $\phi\in H_{\rm per}^{\infty}$ is spectrally stable in the evolution problem (\ref{DGH}) if the spectrum of $J\mathcal{L}$ in $L_{\rm per}^2$ is located on the imaginary axis.
\end{definition}

Let $n(\mathcal{L})$ and $z(\mathcal{L})$ be the numbers of negative and zero eigenvalues of the  linear operator $\mathcal{L}$.

Now we deal with the eigenvalues of the linearized operator $\mathcal{L}$. For this purpose, the following two propositions are needed.

\begin{prop}\label{prop-spectrum1}
\cite{Geyer2021,Neves2009} Let $\mathcal{M}:=-\partial_x^2+Q(x)$ be the Schrodinger operator with the even,
L-periodic, smooth potential $Q$. Assume that $\mathcal{M}w=0$ is satisfied by a linear
combination of two solutions $\varphi_1$ and $\varphi_2$ such that
\be
\varphi_1(x+L)=\varphi_1(x)+\theta\varphi_2(x)
\ee
and
\be
\varphi_2(x+L)=\varphi_2(x)
\ee
with some $\theta\in\mathbb{R}$. Assume that $\varphi_2$ has two zeros on the period of $Q$. Then the zero eigenvalue of
$\mathcal{M}$ in $L^2_{per}$ is simple if $\theta\neq 0$ and double if $\theta=0$. It is the second eigenvalue of $\mathcal{M}$ if $\theta\geq 0$ and the third eigenvalue of $\mathcal{M}$ if $\theta<0.$
\end{prop}

\begin{remark}
From  proposition \ref{prop-spectrum1}, we known

$\bullet\ $ $\theta>0\Rightarrow n(\mathcal{M})=1, z(\mathcal{M})=1$;

$\bullet\ $ $\theta=0\Rightarrow n(\mathcal{M})=1, z(\mathcal{M})=2$;

$\bullet\ $ $\theta<0\Rightarrow n(\mathcal{M})=2, z(\mathcal{M})=1$.

\end{remark}

\begin{prop}\label{prop-spectrum2}
\cite{ Geyer2021,Lopes2003} Let $\mathcal{J}$ be a self-adjoint operator in a Hilbert space $H$ and $\mathcal{S}$ be a bounded invertible operator in $H$. Then $\mathcal{SJS^*}$ and $\mathcal{J}$ have the same inertia, which is the dimension of the negative, null and positive invariant subspace of $H$.
\end{prop}

\begin{lem}\label{L-spectrum}

The linearized operator $\mathcal{L}:L^2_{\rm per}\rightarrow L^2_{\rm per}$ given by (\ref{mathcalL}) admits

$\bullet\ $ two negative eigenvalues and  a simple zero eigenvalue if $\partial_{C_3}\mathfrak{L}>0$;

$\bullet\ $ a simple negative eigenvalue and  a double zero eigenvalue if $\partial_{C_3}\mathfrak{L}=0$;

$\bullet\ $ a simple negative eigenvalue and  a simple zero eigenvalue if $\partial_{C_3}\mathfrak{L}<0$,
where $\mathfrak{L}(C_1,C_2,C_3,b)$ is the period function for the smooth periodic wave $\phi$. Moreover, the rest of the spectrum
of $\mathcal{L}$ in $L^2_{per}$ is strictly positive.
\end{lem}

\textit{Proof.} The method in Ref.\cite{Geyer2021} will be adopted to prove this lemma. Using (\ref{mathcalL}), one can obtain
\be
\begin{split}
\mathcal{L}\phi'&=[-\alpha^2\partial_x(C_1-\phi)\partial_x
+(C_1-C_2-3\phi+\alpha^2\phi'')]\phi'\\
&=-\alpha^2\partial_x[(C_1-\phi)\phi'']+(C_1-C_2
-3\phi+\alpha^2\phi'')\phi'\\
&=-\alpha^2(C_1-\phi)\phi'''+2\alpha^2\phi'\phi''
+(C_1-C_2-3\phi)\phi'\\
&=(-\alpha^2c-\gamma+\alpha^2\phi)\phi'''
+2\alpha^2\phi'\phi''+(c-2\omega-3\phi)\phi'\\
&=0,
\end{split}
\ee
where we have used (\ref{phi1}) in the last step. Therefore, $\phi'\in {\rm{Ker}}(\mathcal{L})\in H^2_{per}$, i.e., $\phi'$ is a solution of $\mathcal{L}\phi=0$. On the other hand, differentiating of the second-order equation (\ref{intg1new}) with respect to $C_3$ yields
\be
\mathcal{L}\partial_{C_3}\phi=0,
\ee
which implies that $\partial_{C_3}\phi$ is the second, linearly indepent solution of $\mathcal{L}\varphi=0$.

Let $y_1,y_2$ be the fundamental set of solutions associated with the equation $\mathcal{L}\varphi=0$ in $H^2_{\rm per}$ satisfying initial conditions
\be
y_1(0)=1,\ y'_1(0)=0,
\ee
and
\be
y_2(0)=0,\ y'_2(0)=1,
\ee
respectively. We set $\phi(0)=\phi(L)=\phi_+$ for the smooth L-periodic solution, where $\phi_+$ is the turning point for the maximum of  $\phi$ in $x$. Hence, we have $\phi'(0)=\phi'(L)=0$ and can define
\be
y_1(x)=\frac{\partial_{C_3}\phi(x)}{\partial_{C_3}\phi_+},\ y_2(x)=\frac{\phi'(x)}{\phi''(0)}.
\ee
Differentiating of the boundary condition $\phi'(L)=0$ yields
\be
y_1'(L)=-\frac{\partial_{C_3}\mathfrak{L}}{\partial_{C_3}\phi_+}\phi''(0).
\ee
Let
\be
\theta=-\frac{\partial_{C_3}
\mathfrak{L}}{\partial_{C_3}\phi_+}\phi''(0),
\ee
then
\be
y_1(x+L)=y_1(x)+\theta y_2(x).
\ee

It should be noticed that $\mathcal{L}$ is not a Schrodinger operator.  To transform the spectral problem $\mathcal{L}v
=\lambda v$ to the spectrum $\mathcal{M}w=\lambda w$ for a Schrodinger operator, we write $\mathcal{L}v=\lambda v$ as the second-order differential equation
\be
\alpha^2(C_1-\phi)v''-\alpha^2\phi'v'+
[(-\alpha^2\phi''+3\phi+C_2-C_1)+\lambda]v=0.
\ee
Let
\be
\begin{split}
p(x)&=\alpha^2(C_1-\phi(x)),\\
q(x)&=-\alpha^2\phi'(x),\\
r(x)&=-\alpha^2\phi''+3\phi+C_2-C_1,
\end{split}
\ee
then we have
\be\label{px}
p(x)v''+q(x)v'+(r(x)+\lambda)v=0.
\ee
Introducing $w(x)$ by
\be
v(x)=w(x)f(x),
\ee
where
\be
f(x)=\sqrt{\frac{C_1-
\phi(0)}{C_1-\phi(x)}},
\ee
then (\ref{px}) becomes
\be\label{px3}
-w''+Q(x)w=\frac{\lambda}{\alpha^2(C_1-\phi(x))}w,
\ee
where
\be
Q(x)=\frac{C_1-C_2-3\phi(x)}{\alpha^2(C_1-\phi(x))}
+\frac{\phi''}{2(C_1-\phi(x))}
-\frac{1}{4}\bigg(\frac{\phi'(x)}{C_1-\phi(x)}\bigg)^2.
\ee

Let
\be
w=\sqrt{\alpha^2(C_1-\phi)}\hat{w},
\ee
then the spectrum problem (\ref{px3}) is equivalent to the spectral problem
\be
\mathcal{SMS}\hat{w}=\lambda\hat{w},
\ee
where $\mathcal{M}:=-\partial_x^2+Q(x)$ is self-adjoint in $L^2_{\rm per}$ and $\mathcal{S}=\sqrt{\alpha^2(C_1-\phi)}$ is a bounded and invertible multiplication operator in $L^2_{\rm per}.$

By proposition \ref{prop-spectrum2}, the number of negative and zero eigenvalues of the spectral problem (\ref{px3}) coincide with those of the operator $\mathcal{M}$. Let
\be
\begin{split}
\varphi_1:=\sqrt{\frac{C_1-\phi(x)}{C_1-\phi(0)}}y_1,\\
\varphi_2:=\sqrt{\frac{C_1-\phi(x)}{C_1-\phi(0)}}y_2.
\end{split}
\ee
One can verify
\be
\begin{split}
\mathcal{M}\varphi_1=0,\ \varphi_1(0)=1,\ \varphi_1'(0)=0,\\
\mathcal{M}\varphi_2=0,\ \varphi_2(0)=0,\ \varphi_2'(0)=1,\\
\end{split}
\ee
and
\be
\varphi_1(x+L)=\varphi_1+\theta\varphi_2(x).
\ee
Moreover, $\varphi_2$ has two zeros in $\mathbb{T}_L$. By proposition \ref{prop-spectrum1}, noting that
\be
{\rm{sign}}(\theta)=-{\rm{sign}}(\partial_a\mathfrak{L}),
\ee
the assertion of Lemma \ref{L-spectrum} holds. This completes the proof of Lemma \ref{L-spectrum}.

Recall that the mass, momentum and energy functionals are given by
\be
\begin{split}
M(u)&=\int_0^L u\d x,\\
E(u)&=\frac{1}{2}\int_0^L(u^2+\alpha^2u_x^2)\d x,\\
F(u)&=\frac{1}{2}\int_0^L(u^3+\alpha^2uu_x^2+2\omega u^2-\gamma u_x^2)\d x.
\end{split}
\ee

Linearization of the mass and momentum functions at the travelling wave with profile
$\phi$ by using the co-periodic perturbation with the profile $w$ yields the constrained subspace of $L^2_{\rm per}$ of the form
\be
X_0\equiv\{w\in L^2_{\rm per}: \langle 1, w \rangle=0,
\langle \phi-\alpha^2\phi'', w \rangle=0 \}.
\ee

The following lemma shows the two constraints
\be
\langle 1, w \rangle=0,
\langle \phi-\alpha^2\phi'', w \rangle=0
\ee
are invariant in the time evolution of the linearized equation $w_t=-J\mathcal{L}w.$

\begin{lem}\label{lem7}
Let $w_0\in H^2_{\rm per}\cap X_0$. If $w\in C^0(\mathbb{R}, H^2_{\rm per}) \cup C^1(\mathbb{R}, H^1_{\rm per})$ is a solution to the linearized DGH equation $w_t=-J\mathcal{L}w$ with initial data $w_0$, then
\be
w(t,\cdot)\in H^2_{\rm per}\cap X_0, \ \ \forall t\in\mathbb{R}.
\ee
\end{lem}

\textit{Proof.} Since $J$ is skew-adjoint and $J1=0$, we have
\be
\frac{\d}{\d t}\langle 1, w \rangle =-\langle 1, J\mathcal{L}w \rangle
=\langle J 1, \mathcal{L}w \rangle=0.
\ee
Similarly,
\be
\begin{split}
\frac{\d}{\d t}\langle \phi-\alpha^2\phi'', w \rangle&=
\langle \phi-\alpha^2\phi'', -J\mathcal{L}w \rangle\\
&=\langle (1-\alpha^2\partial_x^2)\phi, (1-\alpha^2\partial_x^2)^{-1}
\partial_x\mathcal{L}w \rangle\\
&=\langle -\phi',\mathcal{L}w \rangle \\
&=-\langle \mathcal{L}\phi',w \rangle \\
&=0,
\end{split}
\ee
where we have used $\mathcal{L}\phi'=0$ in the last step. It follows from the invariance of the two constraints under the time evolution of the linearized equation that if $w_0\in X_0$, then $w(t,\cdot)\in X_0$ for every $t\in\mathbb{R}$. This completes the proof of Lemma \ref{lem7}.

To study the spectral stability of periodic travelling waves, one need to fix the period $L$ and consider the family of $L-$periodic solutions along a curve in the $(C_2,b)$ plane for fixed $C_1,C_2$ with $2C_1+C_2>0$. The next result describes the character of the $L-$periodic travelling wave solutions, and it is just the first part of Theorem \ref{thm4}.

\begin{lem}\label{L-period}

Fix $C_1,C_2$ with $2C_1+C_2>0$ and $L>0$. There exists a $C^1$ mapping $C_3\mapsto b=\mathfrak{B}_L(C_3)$ for $C_3\in (C_3^L)$ with some $C_3^L\in(0,C_{3,critical})$ and a $C^1$ mapping $a\mapsto\phi=\Phi_L(\cdot,a)\in H^{\infty}_{\rm per}$ of smooth L-periodic solutions along the cuver $b=\mathfrak{B}_L(C_3).$

\end{lem}

\textit{Proof.} For every $C_1,C_2$ with $2C_1+C_2>0$ and $L>0$,  from the monotonicity results obtained in Lemmas (\ref{c3=0}) and (\ref{b-minus}), there exist exactly one $L-$periodic solution on the left and right boundaries of the existence domain on the $(C_3,b)$-plane. The left boundary corresponds to $C_3=0$ and the right boundary corresponds to $C_3=C_3^L$, where $C_3^L$ is uniquely determined by the equation
\be
\mathfrak{L}(C_1,C_2,b_-(C_3^L),C_3^L)=L.
\ee
Since  $\mathfrak{L}(C_1,C_2,b_-(C_3^L)$ is smooth in $(C_1,C_2,b_-(C_3^L)$ and it is strictly increasing in $b$ by Theorem \ref{thm2}, the existence of the $C^1$ mapping $C_3\mapsto b=\mathfrak{B}_L(C_3)$ for $C_3\in(0,C_3^L)$ follows by the implicit function theorem. Indeed, along $b=\mathfrak{B}_L(C_3)$,
\be
\partial_{C_3}\mathfrak{L}+\mathfrak{B}'_L(C_3)\partial_b\mathfrak{L}=0.
\ee
Due to the fact $\partial_b\mathfrak{L}>0$, it follows $\mathfrak{B}'_L(C_3)$ is uniquely defined for every $C_3\in(0,C_3^L)$. Moreover, the mapping $C_3\mapsto\phi=\Phi_L(\cdot,C_3)$ is $C^1$ along the curve  $b=\mathfrak{B}_L(C_3)$, since $\phi$ is smooth with respect to its parameters.
This completes the proof of Lemma \ref{L-period}.

Now it is the time to prove theorem \ref{thm4}. To this end, we introduce the following well-known lemma\cite{Geyer2022,Haragus2008}, see the Remark 5.5 in Ref\cite{Geyer2022}.

\begin{lem}\label{lem-s}
If $n(\mathcal{L}|_{X_0})=0,\ z(\mathcal{L}|_{X_0})=1$, then the spectrum of $J\mathcal{L}$ in $L^2_{\rm per}$ is located on the imaginary axis.
\end{lem}

%Noting that $
%n(\mathcal{L})=0,\ z(\mathcal{L})=1$  is NOT ture for DGH equation, so we adpot A. Geyer and D. Pelinovsky's idea to study
%$n(\mathcal{L}|_{X_0})$ and $z(\mathcal{L}|_{X_0})$\cite{Geyer2022}.

 The counting formulas for the negative and zero eigenvalues of $\mathcal{L}|_{X_0}$ are given by \cite{Geyer2022,Natali2020,Natali2022}
\be
\begin{split}
n(\mathcal{L}|_{X_0})&=n(\mathcal{L})-n_0-z_0,\\
z(\mathcal{L}|_{X_0})&=z(\mathcal{L})+z_0,
\end{split}
\ee
where $n_0$ and $z_0$ are the numbers of negative and zero eigenvalues (taking count their multiplicities) of the matrix of projections
\be
S\equiv \left(\begin{matrix}
\langle\mathcal{L}^{-1}1, 1  \rangle & \langle\mathcal{L}^{-1}(\phi-\alpha^2\phi''), 1  \rangle\\
\langle\mathcal{L}^{-1}1, \phi-\alpha^2\phi''  \rangle & \langle\mathcal{L}^{-1}(\phi-\alpha^2\phi''), \phi-\alpha^2\phi''  \rangle
\end{matrix}\right).
\ee

Formal differentiation of the second-order equation (\ref{intg1}) in $b$ yields
\be
\mathcal{L}(\partial_b\phi)=1,
\ee
which implies
\be
\mathcal{L}^{-1}1=\partial_b\phi.
\ee
Formal differentiation of the second-order equation (\ref{intg1}) in $c$ yields
\be
\mathcal{L}(\partial_c\phi)=-(\phi-\alpha^2\phi''),
\ee
which implies
\be
\mathcal{L}^{-1}(\phi-\alpha^2\phi'')=-\partial_c\phi.
\ee
Therefore the matrix of projections becomes
\be
S\equiv \left(\begin{matrix}
\langle\partial_b\phi, 1  \rangle & \langle-\partial_c\phi, 1  \rangle\\
\langle\partial_b\phi, \phi-\alpha^2\phi''  \rangle & \langle-\partial_c\phi, \phi-\alpha^2\phi''  \rangle
\end{matrix}\right).
\ee
Noting that
\be
\begin{split}
\langle\partial_b\phi, 1  \rangle &=\int_0^L\partial_b\phi\d x=\partial_b \mathcal{M}_L,\\
\langle\partial_b\phi, \phi-\alpha^2\phi''  \rangle&=
\int_0^L(\partial_b\phi)(\phi-\alpha\phi'')\d x\\
&=\int_0^L\phi\partial_b\phi+\alpha^2\phi'\partial_b\phi'\d x\\
&=\partial_b\int_0^L\frac{1}{2}(\phi^2+\alpha^2\phi'^2)\d x\\
&=\partial_b\mathcal{E}_L,
\end{split}
\ee
the matrix $S$ is simply
\be
S\equiv \left(\begin{matrix}
\partial_b\mathcal{M}_L & -\partial_c\mathcal{M}_L\\
\partial_b\mathcal{E}_L & -\partial_c\mathcal{E}_L
\end{matrix}\right).
\ee
Due to the  scaling transformation
\be
\phi(x)=c\hat{\phi}(x), b=c^2\hat{b}, a=c^3\hat{a}, \gamma=c\hat{\gamma}, \omega=c\hat{\omega},
\ee
we have
\be
\mathcal{M}_L(b)=c\hat{\mathcal{M}}_L(\hat{b}),\ \mathcal{E}_L(b)=c^2\hat{\mathcal{E}}_L(\hat{b})
\ee
where the hat quantities are $c-$independent. In terms of the hat quantitied,  the determinant of $S$ reads
\be
\begin{split}
\det S&=c^2(\partial_b\hat{\mathcal{E}}_L)\hat{\mathcal{M}}_L-2c^2
(\partial_b\hat{\mathcal{M}}_L)\hat{\mathcal{E}}_L\\
&=(\partial_{\beta}\hat{\mathcal{E}}_L)\hat{\mathcal{M}}_L-
(\partial_{\beta}\hat{\mathcal{M}}_L)\hat{\mathcal{E}}_L\\
&=\hat{\mathcal{M}}_L^3(\beta)\frac{\d}{\d\beta}\bigg(\frac{\hat{\mathcal{E}}_L(\beta)}{\hat{\mathcal{M}}_L^2(\beta)}\bigg).
\end{split}
\ee

If $\hat{\mathcal{M}}_L^3(\beta)\frac{\d}{\d\beta}\bigg(
\frac{\hat{\mathcal{E}}_L(\beta)}{\hat{\mathcal{M}}_L^2(
\beta)}\bigg)<0,\partial_{C_3}\mathfrak{L}<0$,  then $n_0=1,z_0=0$ and $n(\mathcal{L})=1,z(\mathcal{L})=1$ by Lemma \ref{L-spectrum}.  It follows that
\be
n(\mathcal{L}|_{X_0})=n(\mathcal{L})-n_0-z_0=0,\ z(\mathcal{L}|_{X_0})=z(\mathcal{L})+z_0=1.
\ee
Therefore, by lemma \ref{lem-s}, the spectrum of $J\mathcal{L}$ in $L^2_{per}$ is located on the imaginary axis, which implies that the $L-$period wave is spectrally stable. This complete the proof of Theorem \ref{thm4}.

In the rest of this section, we do a little discussion about the orbital stability. Here we follow the approach in Ref.\cite{Alves2019}, where the following result is given in Proposition 3.8 and Theorem 4.2.

\begin{prop}\label{p4}
Let $V(u)$ be a conserved quantity in the time evolution of the Hamiltonian system
\be
\begin{split}
\frac{\d u}{\d t}=&J\frac{\delta F}{\delta u},\ \ J=-(1-\alpha^2\partial_x^2)^{-1}\partial_x,\\
\frac{\delta F}{\delta u}=&\frac{3}{2}u^2-\frac{1}{2}\alpha^2u_x^2-\alpha^2uu_{xx}+2\omega u+\gamma u_{xx}.
\end{split}
\ee
Assume that the linearized operator $\mathcal{L}$ at the periodic travelling wave with profile
$\phi$ admits a simple negative and a simple zero eigenvalue with
$\rm{ker}(\mathcal{L})=span(\phi')$ satisfying
$\langle V'(\phi),\phi'\rangle=0.$ Assume that there exists
$Y\in H^2_{per}$ such that $\langle \mathcal{L}Y,v\rangle=0$ for every $v\in L^2_{per}$ such that $\langle V'(\phi),v\rangle=0$. If $\langle \mathcal{L}Y,Y\rangle<0$, then the periodic travelling waves is orbitally stable.
\end{prop}

\begin{thm} \label{thm-orbital}
For fixed $C_1,C_2$ which satisfies $2C_1+C_2>0$, $b\leq 0$, $\frac{C_2}{2}[2b+\frac{C_2}{2}(C_1-C_2)]<0$ and $M(\phi)>0$, the smooth periodic wave with profile $\phi\in H^{\infty}_{per}$ of DHG equation is orbitally stable in $H^1_{per}$ when $\partial_{C_3}\mathfrak{L}<0.$

\end{thm}

proof.  Define the conserved quantity $V(u)$ as
\be
V(u):=[2b+\frac{C_2}{2}(C_1-C_2)]M(u)-(C_1+\frac{C_2}{2})E(u),
\ee
then
\be
V'(u)=2b+\frac{C_2}{2}(C_1-C_2)-(C_1+\frac{C_2}{2})(u-\alpha^2u_{xx}).
\ee

Recall that
\be
\mathcal{L}=-\alpha^2\partial_x(C_1-\phi)
\partial_x+(C_1-C_2-3\phi+\alpha^2\phi''),
\ee
then
\be
\begin{split}
\mathcal{L}\phi=&-\alpha^2\partial_x(C_1-\phi)
\partial_x\phi+(C_1-C_2-3\phi+\alpha^2\phi'')\phi\\
=&\alpha^2\phi'^2-(C_1-\phi)\alpha^2\phi''+(C_1
-C_2-3\phi+\alpha^2\phi'')\phi\\
=&2b-2\alpha^2(\phi-C_1)\phi''
-(2C_1-2C_2-3\phi)\phi\\
&-(C_1-\phi)\alpha^2\phi''+(C_1
-C_2-3\phi+\alpha^2\phi'')\phi,\\
=&2b+C_1\alpha^2\phi''-(C_1-C_2)\phi.
\end{split}
\ee

Taking $Y=\phi+\frac{C_2}{2}$, then
\be
\mathcal{L}Y=
2b+\frac{C_2}{2}(C_1-C_2)-(C_1+\frac{C_2}{2})\phi
+(C_1+\frac{C_2}{2})\alpha\phi''.
\ee
Therefore $\langle\mathcal{L}Y,v\rangle=0$ if $\langle V'(\phi),v\rangle=0$. Furthermore
\be
\begin{split}
\langle\mathcal{L}Y,Y\rangle=&\langle
2b+\frac{C_2}{2}(C_1-C_2)-(C_1+\frac{C_2}{2})(\phi-\alpha^2\phi'' ,\phi+\frac{C_2}{2}\rangle\\
=&[2b-\frac{3}{4}C_2^2]M(\phi)+
\frac{C_2}{2}[2b+\frac{C_2}{2}(C_1-C_2)]L
-(C_1+\frac{C_2}{2})E(\phi),
\end{split}
\ee
which implies the periodic waves are stable in the time evolution by Proposition \ref{p4}. This complete the proof of Theorem \ref{thm-orbital}.

%\section{Conclusions and Discussions}\label{S4}

%\appendix

%\section{Proof of Lemma \ref{lem2}}\label{app1}

%\section{Proof of Lemma \ref{lem3}}\label{app2}

%\section{Proof of Lemma \ref{lem4}}\label{app3}

\section*{Acknowledgments}
  X. He is supported by the Natural Science Foundation of Hunan Province (No.2023JJ30179).
A. Chen is supported by the National Natural Science Foundation of China (No.11971163).

%%%%%%%%%%%%%%%%%%%%%%%%%%%%%%%%%%%%%%%%%%%%%%%%%%%%%%%%%%%%%%%%%%%%%%%%%

%%%%%%%%%%%%%%%%%%%%%%%%%%%%%%%%%%%%%%%%%%%%%%%%%%%%%%%%%%%%%%%%%%%%%%%%%
\end{document}